\documentclass[12pt]{amsart}
\usepackage{amsmath,amssymb}

\newtheorem{thm}{Theorem}
\newtheorem{prop}[thm]{Proposition}
\newtheorem{lem}[thm]{Lemma}

\newtheorem{claim}[thm]{Claim}

\newtheorem{rem}[thm]{Remark}

\title{Computation of some leafwise cohomology ring}

\author[Shota MORI]{Shota MORI}
\address[Shota MORI]{Graduate School of Mathematics, Nagoya University, Furocho, Chikusaku, Nagoya, 464-8602, Japan}
\email{mori.shouta.g1@s.mail.nagoya-u.ac.jp}

\date{}


\begin{document}

\maketitle

\begin{abstract}
Let $G$ be the group $SL(2,\mathbb{R})$, $P\subset G$ be the parabolic subgroup of upper triangular matrices and $\Gamma\subset G$ be a cocompact lattice. A right action of $P$ on $\Gamma\backslash G$ defines an orbit foliation $\mathcal{F}_P$.  We compute the leafwise cohomology ring $H^*(\mathcal{F}_P)$ by exploiting non-abelian harmonic analysis on $G$.
\end{abstract}

\section{Introduction and main result}

In the rigidity theory of $C^\infty$ foliations, the leafwise cohomology is a fundamental tool (see Asaoka's survey\cite{asasurvey}). Specifically, in many cases the 1-dimensional leafwise cohomology group plays a crucial role to analize the properties of $C^{\infty}$ foliations. For example, set $G=SL(2,\mathbb{R})$. Let $P\subset G$ be the subgroup of all upper triangular matrices, $\mathfrak{p}$ be the Lie algebra of $P$ and $\Gamma\subset G$ be a cocompact lattice. Also set $M_\Gamma =\Gamma\backslash G$. Let $\mathcal{F}_P$ be the orbit foliation induced from the natural action of $P$ on $M_\Gamma$ and $(\Omega^{*}(\mathcal{F}_P), d_{\mathcal{F}_P})$ be the leafwise complex.  Let $Z^{*}(\mathcal{F}_P)$ be the space of cocycles and $B^{*}(\mathcal{F}_P)$ be the space of coboundaries. Also $H^*(\mathcal{F}_P)$ denote its cohomology group. Then, Matsumoto-Mitsumatsu\cite[Theorem 1]{matmit} proved that
\begin{align}\label{matsumitsuresult}
H^1(\mathcal{F}_P)\cong H^1_{\mathrm{Lie}}(\mathfrak{p}) \oplus H^1_{\mathrm{dR}}(M_\Gamma).
\end{align}

On the other hand, not many examples of higher leafwise cohomology groups are known, except for a linear foliation on a torus\cite{arrsanT}. Recently, Maruhashi-Tsutaya\cite[Theorem 66]{marutsuta} provided a new example. They proved that
\begin{align}\label{marutsutaresult}
H^2(\mathcal{F}_P)\cong H^2_{\mathrm{dR}}(M_\Gamma).
\end{align}
They identified the total leafwise cohomology group $H^*(\mathcal{F}_P)$ as a linear space; however the ring structure remained unknown.

In this paper, we determine the ring structure of $H^*(\mathcal{F}_P)$ by exploiting non-abelian harmonic analysis, which method is totally independent of theirs. Let $g$ be the multiplicity of the irreducible unitary representation $U^{-1}$ whose lowest weight is $1$ in $L^2(M_{\Gamma})$ with the $G$-invariant measure. The finiteness is ensured by a duality theorem\cite[Theorem 1.4.2]{gelpy}. This theorem says that the multiplicity of $U^{1}$ whose highest weight is $-1$ is also $g$. Then, our theorem below holds.
\begin{thm}\label{MyThm2v3}
There exist 1-cocycles $x, y_1, ... , y_{2g}$ such that
\begin{align}
H^{*}(\mathcal{F}_P) \cong \bigwedge \lbrack x, y_{1}, ... , y_{2g} \rbrack \big/ \bigl( \{y_{i} \wedge y_{j} \}_{1 \leq i,j \leq 2g} \bigr).
\end{align}
\end{thm}
In the proof of Theorem \ref{MyThm2v3}, we also prove following: let $C_{\Gamma\,k}=O(2^{5k})$ be the constant defined in (\ref{our_constant_C}) for each $k\in\mathbb{N}$. Then, for each $\eta\in B^*(\mathcal{F}_P)$, there exists $\xi\in \Omega^*(\mathcal{F}_P)$ with $\eta=d_{\mathcal{F}_P}\xi$ such that
\begin{align}
||\xi||_{k}^2 \leq C_{\Gamma \,k+3}||\eta||_{k+3}^2
\end{align}
for each $k\in\mathbb{N}$, where $|| \cdot ||_k$ is a $L^2$-Sobolev norm of $k$-th order. This satisfies a tame estimate introduced in \cite[Definition II.2.1.1.]{hamilton82}. Thus our method is expected to help the further study of foliations using tameness. 

The following describes the flow of this paper. First, in Section \ref{Prel}, we provide tools for discussion. Notations of representation theory, $L^2$-Sobolev norm, and the constants $C_{\Gamma\, k}$ are provided here. Second, in Section  \ref{compsecondcocycles} and \ref{compfirstcocycles}, we compute all explicit generators of $H^*(\mathcal{F}_P)$ in terms of representation theory. We also make a explicit estimations of $L^2$-Soborev norms. Third, we determine the ring structure by computing eigenvalues in section \ref{determinering}.

Finally, the author acknowledges helpful comments from Hirokazu Maruhashi, Shuhei Maruyama, and Hitoshi Moriyoshi. Especially, the author thanks Shuhei Maruyama for pointing out Massey's paper \cite{massey_ring}.

\section{Preliminary}\label{Prel}
We summarize computation formulae.

\subsection{From representation theory}
Let $\hat{G}$ be the unitary dual of $G$. It is sufficient for us to compute $d_{\mathcal{F}_P}$ on each $\pi\in\hat{G}$ by using the differential representation. Indeed, $L^2(M_{\Gamma})$ decomposes into a countable sum of irreducible unitary representations\cite[Theorem 1.2.3]{gelpy}. Set $\mathfrak{g}=\mathfrak{sl}(2,\mathbb{R})$ and take elements
\begin{align}\label{lieelements}
X_0=
\begin{pmatrix}
0 & \frac{1}{2}\\
-\frac{1}{2} & 0
\end{pmatrix}
, X_1=
\begin{pmatrix}
\frac{1}{2} & 0\\
0 & -\frac{1}{2}
\end{pmatrix}
, X_2=
\begin{pmatrix}
0 & 1\\
0 & 0
\end{pmatrix}
\end{align}
from $\mathfrak{g}$. When we regard $X_0, X_1, X_2$ as vector fields on $M_\Gamma$, let $\omega_0, \omega_1, \omega_2\in\Omega^1(M_\Gamma)$ be the dual forms of them. We put
\begin{align}
Y=-X_0+X_2=
\begin{pmatrix}
0 & \frac{1}{2} \\
\frac{1}{2} & 0
\end{pmatrix}
.
\end{align}
For $\pi\in\hat{G}$, let $\pi'$ be the derivative of $\pi$, set
\begin{align}
H=-i\pi'(X_0),\,E=\pi'(X_1)+i\pi'(Y),\,F=-\pi'(X_1)+i\pi'(Y).
\end{align}
When we regard $\pi'$ as the representation of a complex Lie algebra $\mathfrak{sl}(2,\mathbb{C})$, $H$, $E$ and $F$ are given by
\begin{align}
H=\pi'
\begin{pmatrix}
0 & -\frac{i}{2}\\
\frac{i}{2} & 0
\end{pmatrix}
,\,\,E=\pi'
\begin{pmatrix}
\frac{1}{2} & \frac{i}{2}\\
\frac{i}{2} & -\frac{1}{2}
\end{pmatrix}
,\,\,F=\pi'
\begin{pmatrix}
-\frac{1}{2} & \frac{i}{2}\\
\frac{i}{2} & \frac{1}{2}
\end{pmatrix}
.
\end{align}
We set
\begin{align}
C=\pi'(X_0)^2 -\pi'(X_1)^2 -\pi'(Y)^2
\end{align}
which is called the Casimir element. We state the structure theory of $\hat{G}$. Here, $\mathbb{N}$ denotes the set of all non-negative integers and $\frac{1}{2}\mathbb{Z}$ denotes the set of all integers and half-integers.
\begin{thm}\label{SL2Rirruni}(see \cite[Proposition 6.13, Theorem6.2, 6.4, 6.5]{mits})\\
For $\pi\in\hat{G}$, let $\mathbb{M}\subset\frac{1}{2}\mathbb{Z}$ and $q\in\mathbb{R}$ be given in the table below. Then, there exists an orthonormal basis $\{\phi_m\}_{m\in\mathbb{M}}$ of $\pi$ such that
\begin{align}
C\phi_m=q\phi_m,
\end{align}
\begin{align}
H\phi_m=m\phi_m,
\end{align}
\begin{align}
E\phi_m=\sqrt{q+m(m+1)}\phi_{m+1},
\end{align}
\begin{align}
F\phi_m=\sqrt{q+m(m-1)}\phi_{m-1}
\end{align}
for each $m\in\mathbb{M}$.
\begin{table}[htb]
\begin{tabular}{|c|c|c|c|}
\hline
$\pi$ & $\mathbb{M}$ & $q$ & conditions \\
\hline
$V^{0,\frac{1}{2}+i\nu}$ & $\mathbb{Z}$ & $\frac{1}{4}+\nu^2$ & $\nu\geq 0$\\
$V^{\frac{1}{2},\frac{1}{2}+i\nu}$ & $\frac{1}{2}+\mathbb{Z}$ & $\frac{1}{4}+\nu^2$ & $\nu > 0$ \\
$V^{\frac{1}{2},\frac{1}{2}}_+$ & $\frac{1}{2}+\mathbb{N}$ & $\frac{1}{4}$ & -\\
$V^{\frac{1}{2},\frac{1}{2}}_-$ & $-\frac{1}{2}-\mathbb{N}$ & $\frac{1}{4}$ & -\\
$U^{n}$ & $-n-\mathbb{N}$ & $n(1-n)$ & $n\in\frac{1}{2}\mathbb{Z}$ with $n\geq1$\\
$U^{-n}$ & $n+\mathbb{N}$ & $n(1-n)$ & $n\in\frac{1}{2}\mathbb{Z}$ with $n\geq1$\\
$V^{\sigma}$ & $\mathbb{Z}$ & $\sigma(1-\sigma)$ & $\frac{1}{2}<\sigma<1$ \\
$I$ & $\{0\}$ & $0$ & -\\
\hline
\end{tabular}
\end{table}
\end{thm}
The notations $\mathbb{M}$ and $q$ are sometimes denoted as $\mathbb{M}_\pi$ and $q_\pi$, respectively. By applying Theorem \ref{SL2Rirruni}, we write down $d_{\mathcal{F}_P}:\Omega^0(\mathcal{F}_P)\to\Omega^1(\mathcal{F}_P)$. Let $h$ be a $C^\infty$ vector of some $\pi\in\hat{G}$. We have the Fourier expansion 
\begin{align}
h=\sum_{m\in\mathbb{M}}h_{m}\phi_m.
\end{align}
Set $f_1\omega_1+f_2\omega_2=d_{\mathcal{F}_P}h$. Then, Fourier coefficients of $f_1$ and $f_2$ are given by
\begin{align}\label{compdFP01_f1}
f_{1\,m}=\frac{h_{m-1}}{2}\alpha_{q\,m-1} - \frac{h_{m+1}}{2}\beta_{q\,m+1}
\end{align}
and
\begin{align}\label{compdFP01_f2}
f_{2\,m}=-\frac{ih_{m-1}}{2}\alpha_{q\,m-1}+imh_{m}-\frac{ih_{m+1}}{2}\beta_{q\,m+1},
\end{align}
where
\begin{align}
\alpha_{q\,m}=\sqrt{q+m(m+1)},
\end{align}
\begin{align}
\beta_{q\,m}=\sqrt{q+m(m-1)}.
\end{align}

Next, we write down $d_{\mathcal{F}_P}:\Omega^1(\mathcal{F}_P)\to\Omega^2(\mathcal{F}_P)$. Let $f_1, f_2$ be $C^\infty$ vectors of some $\pi\in\hat{G}$. Set $g\,\omega_1\wedge\omega_2=d_{\mathcal{F}_P}(f_1\omega_1+f_2\omega_2)$. Then, $g$'s Fourier coefficients are given by
\begin{align}\label{compdFP12}
g_m=\frac{if_{1\,m-1}+f_{2\,m-1}}{2}\alpha_{q\,m-1}-(imf_{1\,m}+f_{2\,m})+\frac{if_{1\,m+1}-f_{2\,m+1}}{2}\beta_{q\,m+1}.
\end{align}

For convenience, we replace $h$ by $-4h$, $f_{1\,m}$ by $2f_{1\,m}$, $f_{2\,m}$ by $2if_{2\,m}$, and $g$ by $ig$. Then we always assume
\begin{align}\label{compdFP01_f1_changed}
f_{1\,m}=-h_{m-1}\alpha_{q\,m-1} + h_{m+1}\beta_{q\,m+1},
\end{align}
\begin{align}\label{compdFP01_f2_changed}
f_{2\,m}=h_{m-1}\alpha_{q\,m-1}-2mh_{m}+h_{m+1}\beta_{q\,m+1},
\end{align}
or
\begin{align}\label{compdFP12_changed}
g_m=(f_{1\,m-1}+f_{2\,m-1})\alpha_{q\,m-1}-2(mf_{1\,m}+f_{2\,m})+(f_{1\,m+1}-f_{2\,m+1})\beta_{q\,m+1}.
\end{align}

\subsection{From Sobolev space theory}
In this paper, we construct formal functions by using (\ref{compdFP01_f1_changed}), (\ref{compdFP01_f2_changed}), and (\ref{compdFP12_changed}). To ensure their smoothness, we use $L^2$-Sobolev norms. In general, let $M$ be a compact Riemannian manifold and $(\lambda_s)_{s=0}^{\infty}$ be the sequence consisting of eigenvalues of the Laplace-Beltrami operator:
\begin{align}
0=\lambda_0 \leq \lambda_1 \leq ... \leq \lambda_s \leq ... \to \infty.
\end{align}
Then, for each $k\in\mathbb{N}$, the $L^2$-Sobolev norm of $k$-th order with respect to the Bessel potential is given by
\begin{align}
||f||_k^2=\sum_{s=0}^{\infty}(1+\lambda_s)^k|f_s|^2,
\end{align}
where $f\in C^\infty(M)$ and $(f_s)_{s=0}^{\infty}$ is the Fourier coefficients of $f$. (For example, see \cite[Definition 4.1]{strichartz}).

We apply this fact to our case. We define the Riemannian metric on $M_\Gamma$ whose orthogonal frame is $\{X_0, X_1, Y\}$. Then the Laplace-Beltrami operator is $\Delta=-X_0^2-X_1^2-Y^2$. (See \cite[Theorem 1]{ura}. Since $G$ is connected unimodular, we have $\mathrm{Trace}(\mathrm{ad}(\cdot))=0$.) We transform with $\Delta=C-2X_0^2$. Then we get
\begin{align}
\Delta\phi_m = (q+2m^2)\phi_m.
\end{align}
Thus a $L^2$-Sobolev norm of $k$-th order is given by
\begin{align}\label{our_Sobolev}
||f||_k^2=\sum_{\pi\subset L^2(M_\Gamma), \,m\in\mathbb{M}_\pi}(1+q+2m^2)^k|f_{\pi\, m}|^2
\end{align}
for each $k\in\mathbb{N}$ and $f\in C^\infty(M_\Gamma)$.

We estimate the $L^2$-Sobolev norms (\ref{our_Sobolev}) in Sections \ref{compsecondcocycles} and \ref{compfirstcocycles}. We provide constants for this purpose. Observe the behavior of number $q$:
\begin{lem}\label{q_dist}
The sequence $(q_\pi)_{\pi\subset L^2(M_\Gamma)}$ has no accumulation points, where $q_\pi$ is the number of $\pi$.
\end{lem}
\begin{proof}
The non-existence of accumulation points is derived from the fact that eigenvalues of $\Delta$ diverge.
\end{proof}
Set
\begin{align}\label{q_Gamma}
q_\Gamma = \inf \left\{ |q_\pi| \, \middle| \, \pi\subset L^2(M_\Gamma), \,\, \pi \neq I, U^{-1}, U^{1} \right\} >0 \,\,\,\text{(Lemma \ref{q_dist})}.
\end{align}
If $q_\Gamma>1$, then we replace $q_\Gamma = 1$.

For each $k\in\mathbb{N}$, set
\begin{align}\label{our_constant_C}
C_{\Gamma\, k}= \Bigl( \frac{3!}{q_\Gamma^2}\Bigr)^2 2^{5k+8}.
\end{align}

The estimation is performed using this constant $C_{\Gamma\, k}$ under replacing $h$, $f_1$, $f_2$, and $g$. When a special cocycle $\eta\in Z^*(\mathcal{F}_P)$ is given, we construct $\xi\in\Omega^*(\mathcal{F}_P)$ which satisfies $\eta=d_{\mathcal{F}_P}\xi$ formally and prove
\begin{align}
||\xi||_{k}^2 \leq C_{\Gamma \,k+3}||\eta||_{k+3}^2
\end{align}
for each $k\in\mathbb{N}$.

\section{Computing second cocycles}\label{compsecondcocycles}
To prove the lemmata below, we solve a linear equation for all Fourier coefficients on each $\pi\in\hat{G}$. The symbol $|_\pi$ means ``restricted to $\pi$".

\subsection{Trivial representation}\label{trivcomp_2cocycle}
We get the following lemma directly.
\begin{lem}\label{trivcomplem}
$H^2(\mathcal{F}_P)|_{I}=\{0\}$.
\end{lem}

\subsection{Corresponding to the lowest weight 1}\label{1comp_2cocycle}
We characterize the coboundary space.
\begin{lem}\label{1complem}
The following holds.
\begin{align}
B^2(\mathcal{F}_P)|_{U^{-1}}=\left\{g\,\omega_1\wedge\omega_2\middle|\sum_{m=1}^{\infty}\sqrt{m}g_m=0\right\}.
\end{align}
Specially, $Z^2(\mathcal{F}_P)|_{U^{-1}}$ is spanned by $\phi_1\,\omega_1\wedge\omega_2$ and $B^2(\mathcal{F}_P)|_{U^{-1}}$.
\end{lem}
\begin{proof}
Proving that the right hand side contains $B^2(\mathcal{F}_P)|_{U^{-1}}$ is easy. To prove opposite, let $N$ be a positive integer. We put
\begin{align}
f_{1\,m}^{(N+1)}=\frac{N+1-m}{\sqrt{m(N+1)}}\,\,\,(1\leq m\leq N),
\end{align}
\begin{align}
f_{1\,N+1}^{(N+1)}=0,
\end{align}
\begin{align}
f_1^{(N+1)}=\sum_{m=1}^{N}f_{1\,m}^{(N+1)}\phi_m.
\end{align}
We have
\begin{align}\label{U-1_N+1}
d_{\mathcal{F}_P}(f_1^{(N+1)}\omega_1)=(-\sqrt{N+1}\phi_1+\phi_{N+1})\, \omega_1\wedge\omega_2
\end{align}
by (\ref{compdFP12_changed}). To check it, put $g^{(N+1)}\, \omega_1 \wedge \omega_2 = d_{\mathcal{F}_P}(f_1^{(N+1)}\omega_1)$. For each $2\leq m \leq N$,
\begin{align}
&\, g^{(N+1)}_m  \\
&= f^{(N+1)}_{1\, m-1}\alpha_{q\, m-1}-2mf^{(N+1)}_{1\,m} + f^{(N+1)}_{1 \, m+1}\beta_{q \, m+1} \nonumber \\
&= \frac{N+1-(m-1)}{\sqrt{(m-1)(N+1)}}\sqrt{(m-1)m} -2m\frac{N+1-m}{\sqrt{m(N+1)}} + \frac{N+1-(m+1)}{\sqrt{(m+1)(N+1)}}\sqrt{m(m+1)} \nonumber \\
&= \sqrt{\frac{m}{N+1}} \Bigl( N+1-(m-1)-2(N+1-m)+N+1-(m+1) \Bigr) \nonumber \\
&= 0. \nonumber
\end{align}
Next,
\begin{align}
g^{(N+1)}_1 &= -2\cdot 1 \cdot f^{(N+1)}_{1\,1} + f^{(N+1)}_{1 \, 2}\beta_{q \, 2} \\
&= -2\frac{N+1-1}{\sqrt{1\cdot (N+1)}} + \frac{N+1-2}{\sqrt{2\cdot (N+1)}}\sqrt{1\cdot 2} \nonumber \\
&= -\sqrt{N+1}. \nonumber
\end{align}
Finally,
\begin{align}
g^{(N+1)}_{N+1} &= f^{(N+1)}_{1\, N}\alpha_{q\, N}-2(N+1)f^{(N+1)}_{1\,N+1} \\
&= \frac{N+1-N}{\sqrt{N(N+1)}}\sqrt{N(N+1)} - 0 \nonumber \\
&= 1. \nonumber
\end{align}
Thus the formula (\ref{U-1_N+1}) is valid. We put $\xi^{(N+1)}=f_1^{(N+1)}\omega_1$.

Then let $\eta = g\,\omega_1\wedge\omega_2$ be an element from the right hand side. Put
\begin{align}\label{xi_U-1}
\xi = \sum_{N=1}^{\infty} g_{N+1}\,\xi^{(N+1)}.
\end{align}
We obtain $\eta = d\xi$ formally. This is determined to be smooth 2-coboundary after the Sobolev estimation (Lemma \ref{Sobolev_estimation_U-1}).
\end{proof}
\begin{lem}\label{Sobolev_estimation_U-1}
Let $\eta$ be an element from the right hand side in Lemma \ref{1complem} and put $\xi$ as (\ref{xi_U-1}). Then, for each $k\in\mathbb{N}$, we have
\begin{align}
||\xi||_k^2 \leq ||\eta||_{k+3}^2.
\end{align} 
\end{lem}
\begin{proof}
Fix $k\in\mathbb{N}$. It is enough to prove
\begin{align}
\sum_{m=1}^{\infty} (1+2m^2)^k \left| \sum_{N=m}^{\infty} g_{N+1}\frac{N+1-m}{\sqrt{m(N+1)}} \right|^2 \leq \sum_{N=1}^{\infty} (1+2(N+1)^2)^{k+3} |g_{N+1}|^2.
\end{align}
First, we estimate each term in the left hand side. Put
\begin{align}
g'_{N+1}=(1+2(N+1)^2)g_{N+1}.
\end{align}
Then
\begin{align}
&\, \left| \sum_{N=m}^{\infty} g_{N+1}\frac{N+1-m}{\sqrt{m(N+1)}} \right|^2  \\
&= \left| \sum_{N=m}^{\infty} \frac{g'_{N+1}}{1+2(N+1)^2} \frac{N+1-m}{\sqrt{m(N+1)}} \right|^2 \nonumber \\
&\leq \left( \sum_{N=m}^{\infty} \frac{1}{(1+2(N+1)^2)^2} \frac{(N+1-m)^2}{m(N+1)} \right) \sum_{N=m}^{\infty} |g'_{N+1}|^2 \,\,\, \text{(Cauchy-Schwartz)} \nonumber \\
&\leq \left( \sum_{N=m}^{\infty} \frac{1}{1+2(N+1)^2} \right) \sum_{N=m}^{\infty} |g'_{N+1}|^2 \nonumber \\
&\leq 1 \cdot \sum_{N=m}^{\infty} (1+2(N+1)^2)^2 |g_{N+1}|^2. \nonumber
\end{align}
Next, we estimate the whole term in the left hand side.
\begin{align}
&\, \sum_{m=1}^{\infty} (1+2m^2)^k \left| \sum_{N=m}^{\infty} g_{N+1}\frac{N+1-m}{\sqrt{m(N+1)}} \right|^2 \\
&\leq \sum_{m=1}^{\infty} (1+2m^2)^k \sum_{N=m}^{\infty} (1+2(N+1)^2)^2 |g_{N+1}|^2 \nonumber \\
&= \sum_{N=1}^{\infty} (1+2(N+1)^2)^2 \left( \sum_{m=1}^{N} (1+2m^2)^k \right) |g_{N+1}|^2 \nonumber \\
&\leq \sum_{N=1}^{\infty} (1+2(N+1)^2)^2 N (1+2N^2)^k |g_{N+1}|^2 \nonumber \\
&\leq \sum_{N=1}^{\infty} (1+2(N+1)^2)^{k+3}|g_{N+1}|^2. \nonumber
\end{align}
\end{proof}

\subsection{Corresponding to the highest weight -1}\label{-1comp_2cocycle}
Similar argument also holds in this case.
\begin{lem}\label{-1complem}
The following holds.
\begin{align}
B^2(\mathcal{F}_P)|_{U^{1}}=\left\{g\,\omega_1\wedge\omega_2\middle|\sum_{m=-\infty}^{-1}(-1)^{-m}\sqrt{-m}g_m=0\right\}.
\end{align}
Specially, $Z^2(\mathcal{F}_P)|_{U^{1}}$ is spanned by $\phi_{-1}\,\omega_1\wedge\omega_2$ and $B^2(\mathcal{F}_P)|_{U^{1}}$.
\end{lem}
\begin{proof}
Proving that the right hand side contains $B^2(\mathcal{F}_P)|_{U^{1}}$ is easy. To prove opposite, let $N$ be a positive integer. We put
\begin{align}
f_{1\,m}^{-(N+1)}=\frac{N+1+m}{\sqrt{-m(N+1)}}\,\,\,(-N\leq m\leq -1),
\end{align}
\begin{align}
f_{1\, -(N+1)}^{-(N+1)} = 0,
\end{align}
\begin{align}
f_1^{-(N+1)}=\sum_{m=-N}^{-1}(-1)^{N-m} f_{1\,m}^{-(N+1)}\phi_m.
\end{align}
Then, we have
\begin{align}\label{U+1_N+1}
d_{\mathcal{F}_P}(f_1^{-(N+1)}\omega_1)=(\phi_{-(N+1)}+(-1)^{N+1}\sqrt{N+1}\phi_{-1})\, \omega_1\wedge\omega_2
\end{align}
by (\ref{compdFP12_changed}). To check it, put $g^{-(N+1)}\, \omega_1 \wedge \omega_2 = d_{\mathcal{F}_P}(f_1^{-(N+1)}\omega_1)$. For each $-N\leq m \leq -2$,
\begin{align}
&\, (-1)^{N-(m-1)}g^{-(N+1)}_m  \\
&=  f^{-(N+1)}_{1\, m-1}\alpha_{q\, m-1} + 2mf^{-(N+1)}_{1\,m} + f^{-(N+1)}_{1 \, m+1}\beta_{q \, m+1} \nonumber \\
&=  \frac{N+1+(m-1)}{\sqrt{|m-1|(N+1)}}\sqrt{(m-1)m} +\frac{2m(N+1+m)}{\sqrt{|m|(N+1)}} + \frac{N+1+(m+1)}{\sqrt{|m+1|(N+1)}}\sqrt{m(m+1)} \nonumber \\
&= \sqrt{\frac{|m|}{N+1}} \Bigl( N+1+(m-1)-2(N+1+m)+N+1+(m+1) \Bigr) \nonumber \\
&= 0. \nonumber
\end{align}
Next,
\begin{align}
g^{-(N+1)}_{-(N+1)} &= -2(-(N+1))(-1)^{N+(N+1)}f^{-(N+1)}_{1\,-(N+1)} + (-1)^{N+N}f^{(N+1)}_{1 \, -N}\beta_{q \, -N} \\
&= 0 + \frac{N+1-N}{\sqrt{N \cdot (N+1)}}\sqrt{(-N-1)(-N)} \nonumber \\
&= 1. \nonumber
\end{align}
Finally,
\begin{align}
g^{-(N+1)}_{-1} &= (-1)^{N+2}f^{-(N+1)}_{1\, -2}\alpha_{q\, -2}-2(-1)(-1)^{N+1}f^{-(N+1)}_{1\,-1} \\
&= (-1)^{N+1}\left( - \frac{N+1-2}{\sqrt{2\cdot (N+1)}}\sqrt{2\cdot 1} + 2 \frac{N+1-1}{\sqrt{1\cdot (N+1)}} \right) \nonumber \\
&= (-1)^{N+1}\sqrt{N+1}. \nonumber
\end{align}
Thus the formula (\ref{U+1_N+1}) is valid. This fact proves the result as in Lemma \ref{1complem}. The Sobolev estimation is the same as Lemma \ref{Sobolev_estimation_U-1}:
\begin{align}
||\xi||_k^2 \leq ||\eta||_{k+3}^2.
\end{align}
\end{proof}

\subsection{The other cases}\label{othercomp_2cocycle}
Set $\pi\neq I, U^{-1}, U^1$ and fix $m_\pi \in\mathbb{M}$. To begin with, we observe a behaver of $d_{\mathcal{F}_P}|_\pi$. Put $f_{1\, m}=f_{2\,m}=f_{1\, m+3}=f_{2\,m+3}=0$ for any $m\in\mathbb{M}$ which satisfies $m\equiv m_\pi \,\,(\mathrm{mod}\,4)$. We consider the linear map
\begin{align}
(f_{1\,m+1},f_{2\,m+1},f_{1\,m+2},f_{2\,m+2})\mapsto(g_{m},g_{m+1},g_{m+2},g_{m+3})
\end{align}
for all $m\equiv m_\pi$ by (\ref{compdFP12_changed}). The coefficient matrix is a block diagonal matrix whose block is a $4 \times 4$ matrix. Each block is represented as
\begin{align}\label{d12_matrix}
\begin{pmatrix}
g_{m} \\
g_{m+1} \\
g_{m+2} \\
g_{m+3}
\end{pmatrix}
=
\begin{pmatrix}
\beta_{q\,m+1} & -\beta_{q\,m+1} & 0 & 0 \\
-2(m+1) & -2 & \beta_{q\,m+2} & -\beta_{q\,m+2} \\
\alpha_{q\,m+1} & \alpha_{q\,m+1} & -2(m+2) & -2 \\
0 & 0 & \alpha_{q\,m+2} & \alpha_{q\,m+2}
\end{pmatrix}
\begin{pmatrix}
f_{1\,m+1} \\
f_{2\,m+1} \\
f_{1\,m+2} \\
f_{2\,m+2}
\end{pmatrix}
.
\end{align}
We denote this $4 \times 4$ matrix as $A_{q\, m}$. Its determinant $\gamma_{q\,m}=\det A_{q\,m}$ is
\begin{align}\label{our_each_gamma}
\gamma_{q\,m}=-4q\sqrt{m^2+m+q}\sqrt{m^2+5m+6+q}.
\end{align}
This value is always non-vanishing when $m,m+3\in\mathbb{M}$ by the lemma below.
\begin{lem}\label{gamma_and_q}
The following holds:
\begin{align}
|\gamma_{q\, m}| \geq 4 \min \{q^2, 1\}.
\end{align}
Especially, we obtain $|\gamma_{q\, m}| \geq 4q_\Gamma^2$.
\end{lem}
\begin{proof}
When $\pi\neq I, U^{-n}, U^{n}$, since $\mathbb{M}\subset\frac{1}{2}\mathbb{Z}$, we have
\begin{align}
|\gamma_{q\,m}| &= 4q\sqrt{\Bigl(m+\frac{1}{2}\Bigr)^2-\frac{1}{4}+q}\sqrt{\Bigl(m+\frac{5}{2}\Bigr)^2-\frac{1}{4}+q} \\
&\geq 4q^2. \nonumber
\end{align}
When $\pi=U^{-n}\,\,(n\geq \frac{3}{2})$, $|\gamma_{q\,m}|$ takes the minimum at $m=n$:
\begin{align}
|\gamma_{q\,m}| &\geq 4|q|\sqrt{n^2+n+n(1-n)}\sqrt{n^2+5n+6+n(1-n)}  \\
&= 4|q|\sqrt{2n}\sqrt{6(n+1)} \nonumber \\
&\geq 4\cdot \frac{3}{2}\cdot\frac{1}{2}\sqrt{3}\sqrt{15} \nonumber \\
&\geq 4 \cdot 1. \nonumber
\end{align}

When $\pi=U^{n}\,\,(n\geq \frac{3}{2})$, $|\gamma_{q\,m}|$ takes the minimum at $m+3=-n$:
\begin{align}
|\gamma_{q\,m}| &\geq 4|q|\sqrt{(n+3)^2-(n+3)+n(1-n)}\sqrt{(n+3)^2-5(n+3)+6+n(1-n)} \\
&= 4|q|\sqrt{6(n+1)}\sqrt{2n} \nonumber \\
&\geq 4 \cdot 1. \nonumber
\end{align}
\end{proof}
Therefore, we can determine the values $f_{1\,m+1}, f_{2\,m+1}, f_{1\,m+2}, f_{2\,m+2}$ which satisfy (\ref{d12_matrix}). 

Then take any $\eta= g\, \omega_1 \wedge \omega_2 \in Z^2(\mathcal{F}_P)|_\pi$. Put $f_{1\, m}=f_{2\,m}=f_{1\, m+3}=f_{2\,m+3}=0$ for any $m\in\mathbb{M}$ which satisfies $m\equiv m_\pi \,\,(\mathrm{mod}\,4)$. We determine $f_{1\,m+1}, f_{2\,m+1}, f_{1\,m+2}, f_{2\,m+2}$ in (\ref{d12_matrix}) and set
\begin{align}
\xi = f_1\, \omega_1 + f_2\, \omega_2.
\end{align}
We get $\eta=d_{\mathcal{F}_P}\xi$ formally. We last the Sobolev estimations (Lemma \ref{Sobolev_estimation_other} and \ref{Sobolev_estimation_U-+n}).
\begin{lem}\label{Sobolev_estimation_other}
When $\pi \neq I, U^{-n}, U^{n}$, we have
\begin{align}
||\xi||_k^2 \leq  \Bigl( \frac{3!}{q_\Gamma^2}\Bigr)^2 2^{5(k+3)+8} ||\eta||_{k+3}^2
\end{align}
for each $k\in\mathbb{N}$, where $q_\Gamma$ is defined in (\ref{q_Gamma}).
\end{lem}
\begin{proof}
Fix $k\in\mathbb{N}$ and $m\equiv m_\pi$. Each entry of $A_{q\, m}$ is bounded from $2\sqrt{q+m^2}$ or $2\sqrt{q+(m+3)^2}$. Once we assume the latter. Here, a cofactor of $A_{q\, m}$ is degree $3$ polynomial of entries of $A_{q\,m}$. Thus any entry of the cofactor matrix of $A_{q\,m}$ is bounded from
\begin{align}
3! 2^3 (q+(m+3)^2)^{\frac{3}{2}}.
\end{align}
Then, for each $l=1,2$ and $m'=m+1, m+2$,
\begin{align}
|f_{l\, m'}| &\leq \frac{1}{|\gamma_{q\, m}|} 3! 2^3 (q+(m+3)^2)^{\frac{3}{2}} \left( \sum_{m''=m}^{m+3} |g_{m''}| \right) \\
&\leq \frac{3!}{4q_\Gamma^2} 2^3 (q+(m+3)^2)^{\frac{3}{2}} \left( 4\sqrt{\sum_{m''=m}^{m+3} |g_{m''}|^2} \right) \nonumber \\
&\leq \frac{3!}{q_\Gamma^2} 2^3 (1+q+2(m+3)^2)^{\frac{3}{2}} \sqrt{\sum_{m''=m}^{m+3} |g_{m''}|^2} \nonumber
\end{align}
from Lemma \ref{gamma_and_q}.
\begin{claim}\label{Sobolev_change_claim_other}
For each $s', s'' \in \{0, 1, 2, 3\}$, we have
\begin{align}
1+q+2(m+s')^2 \leq 2^5 (1+q+2(m+s'')^2).
\end{align}
\end{claim}
We continue by assuming this claim.
\begin{align}
&\, (1+q+2m'^2)^k|f_{l\,m}|^2 \\
&\leq \Bigl( \frac{3!}{q_\Gamma^2} \Bigr)^2 2^6 (1+q+2m'^2)^k (1+q+2(m+3)^2)^{3} \sum_{m''=m}^{m+3} |g_{m''}|^2 \nonumber \\
&\leq \Bigl( \frac{3!}{q_\Gamma^2} \Bigr)^2 2^{5(k+3)+6} \sum_{m''=m}^{m+3} (1+q+2m''^2)^{k+3}|g_{m''}|^2 \,\,\, \text{(Claim \ref{Sobolev_change_claim_other})}. \nonumber
\end{align}
Add together for $l=1,2$ and $m'=m+1, m+2$ , and get
\begin{align}
&\, \sum_{l=1,2, \,\, m'=m+1,m+2}(1+q+2m'^2)^k|f_{l\,m}|^2 \\
&\leq \Bigl( \frac{3!}{q_\Gamma^2} \Bigr)^2 2^{5(k+3)+8} \sum_{m''=m}^{m+3} (1+q+2m''^2)^{k+3}|g_{m''}|^2 \nonumber
\end{align}
Finally, the desired inequality is obtained by adding up for $m\equiv m_\pi$.
\end{proof}
\begin{proof}[Proof of Claim \ref{Sobolev_change_claim_other}]
We find a constant $c>1$ which satisfies
\begin{align}
1+q+2(m+s)'^2 \leq c (1+q+2(m+s'')^2).
\end{align}
By transposition, we obtain
\begin{align}
2c(m+s'')^2 - 2(m+s')^2 + (c-1)(1+q) \geq 0.
\end{align}
Since $\pi \neq I, U^{-n}, U^{n}$, $q$ is positive. Then it is enough to satisfy
\begin{align}
2c(m+s'')^2 - 2(m+s')^2 + c-1 \geq 0.
\end{align}
In the left hand side, we obtain
\begin{align}
2c(m+s'')^2 - 2(m+s')^2 &= 2(c-1)\left( m+ \frac{cs''-s}{c-1} \right)^2 -2c\frac{(s''-s')^2}{c-1} \\
&\geq -\frac{18c}{c-1}. \nonumber
\end{align}
Then it is enough to satisfy
\begin{align}
-\frac{18c}{c-1} + c-1 \geq 0
\end{align}
or
\begin{align}
-18c+(c-1)^2 \geq 0.
\end{align}
Roughly, this is valid for $c\geq 20$. Thus we can set $c=2^5$.
\end{proof}

\begin{lem}\label{Sobolev_estimation_U-+n}
When $\pi= U^{-n}, U^{n}\,\,(n\geq\frac{3}{2})$, we have
\begin{align}
||\xi||_k^2 \leq  \Bigl( \frac{3!}{q_\Gamma^2}\Bigr)^2 2^{5(k+3)+8} ||\eta||_{k+3}^2
\end{align}
for each $k\in\mathbb{N}$.
\end{lem}
\begin{proof}
We prove the case $\pi=U^{-n}$. (The proof for case $\pi=U^{n}$ is similar.) Fix $k\in\mathbb{N}$ and $m\equiv m_\pi$. Each entry of $A_{q\, m}$ is bounded from $2(m+3)$. Then any entry of the cofactor matrix of $A_{q\,m}$ is bounded from
\begin{align}
3! 2^3 (m+3)^3.
\end{align}
Then, for each $l=1,2$ and $m'=m+1, m+2$,
\begin{align}
|f_{l\, m'}| &\leq \frac{1}{|\gamma_{q\, m}|} 3! 2^3 (m+3)^3 \left( \sum_{m''=m}^{m+3} |g_{m''}| \right)  \\
&= \frac{3!}{q_\Gamma^2} 2^3 (m+3)^3 \sqrt{\sum_{m''=m}^{m+3} |g_{m''}|^2} \nonumber
\end{align}
from Lemma \ref{gamma_and_q}.
\begin{claim}\label{Sobolev_change_claim1_U-+n}
The following holds.
\begin{align}
m+3 \leq \sqrt{1+n(1-n)+2(m+3)^2}.
\end{align}
\end{claim}
We continue by assuming this claim.
\begin{align}
|f_{l\, m'}| \leq \frac{3!}{q_\Gamma^2} 2^3 (1+n(1-n)+2(m+3)^2)^{\frac{3}{2}} \sqrt{\sum_{m''=m}^{m+3} |g_{m''}|^2}\,\,\, \text{(Claim \ref{Sobolev_change_claim1_U-+n})}.
\end{align}
\begin{claim}\label{Sobolev_change_claim2_U-+n}
For each $s', s'' \in \{0, 1, 2, 3\}$, we have
\begin{align}
1+n(1-n)+2(m+s')^2 \leq 2^5 (1+n(1-n)+2(m+s'')^2).
\end{align}
\end{claim}
Assume this claim. The rest is the same as after Claim \ref{Sobolev_change_claim_other} of the proof of Lemma \ref{Sobolev_estimation_other}.
\end{proof}
\begin{proof}[Proof of Claim \ref{Sobolev_change_claim1_U-+n}]
Transforming the desired inequality, we obtain
\begin{align}
1+n(1-n)+(m+3)^2 \geq 0.
\end{align}
The left hand side take the minimum if $m=n$. Then
\begin{align}
1+n(1-n)+(n+3)^2 = 7n+10.
\end{align}
This is positive.
\end{proof}
\begin{proof}[Proof of Claim \ref{Sobolev_change_claim2_U-+n}]
We find a constant $c>13$ which satisfies
\begin{align}
1+n(1-n)+2(m+s')^2 \leq c (1+n(1-n)+2(m+s'')^2).
\end{align}
By transposition, we obtain
\begin{align}
2c(m+s'')^2 - 2(m+s')^2 + (c-1)(1+n(1-n)) \geq 0.
\end{align}
Since $c>3$ and $m \geq n \geq \frac{3}{2}$, we have
\begin{align}
2c(m+s'')^2 - 2(m+s')^2 &\geq 2cm^2-2(m+3)^2  \\
&= 2(c-1)\left( m-\frac{3}{c-1} \right)^2 -\frac{18c}{c-1} \nonumber \\
&\geq 2cn^2-2(n+3)^2. \nonumber
\end{align}
Then it is enough to satisfy
\begin{align}
2cn^2 - 2(n+3)^2 + (c-1)(1+n(1-n)) \geq 0.
\end{align}
Since $c>13$ and $n\geq \frac{3}{2}$, the left hand side is estimated by the below:
\begin{align}
&\, 2cn^2 - 2(n+3)^2 + (c-1)(1+n(1-n))  \\
&= (c-1)n^2+(c-13)n+(c-19) \nonumber \\
&\geq (c-1)\Bigl( \frac{3}{2} \Bigr)^2+(c-13)\frac{3}{2}+(c-19) \nonumber \\
&= \frac{19}{4}c - \frac{163}{4}. \nonumber
\end{align}
This is positive under $c>13$. Thus we can set $c=2^5$ roughly.
\end{proof}

Then the following holds.
\begin{lem}\label{othercomplem}
$H^2(\mathcal{F}_P)|_\pi =\{0\}$.
\end{lem}

\subsection{The whole sum}\label{whole_sum_2coboundaries}
For any $\eta \in B^2(\mathcal{F}_P)$ and $\pi\subset L^2(M_\Gamma)$, let $\eta_\pi$ be the $\pi$-component. Then we can get smooth cochain $\xi_\pi \in \Omega^1(\mathcal{F}_P)|_\pi$ such that 
\begin{align}
|| \xi_\pi ||_k^2 \leq C_{\Gamma \, k+3} || \eta_\pi ||_{k+3}^2
\end{align}
for each $k\in\mathbb{N}$. Thus we get $\eta=d_{\mathcal{F}_P}\xi$ and
\begin{align}
|| \xi ||_k^2 \leq C_{\Gamma \, k+3} || \eta ||_{k+3}^2,
\end{align}
where $\xi = \sum_{\pi} \xi_\pi$.

We summarize the discussion so far. We put
\begin{align}\label{1-cocycle_x}
x=\omega_1,
\end{align}
\begin{align}\label{1-cocycle_y_j}
y_j = \phi_1(\omega_1 - \omega_2) \,\, \text{in $j$-th} \,\, U^{-1},
\end{align}
\begin{align}\label{1-cocycle_y_g_j}
y_{g+j} = \phi_{-1}(\omega_1 + \omega_2) \,\, \text{in $j$-th} \,\, U^{1}.
\end{align}
These are 1-cocycles. Then, we got the following.
\begin{prop}\label{2-cohomology}
The set $\{x\wedge y_1,\, ...\, ,\, x\wedge y_{2g}\}$ is basis for $H^2(\mathcal{F}_P)$, where the number $g$ is the multiplicity of $U^{-1}$ and $U^{1}$.
\end{prop}
\begin{rem}
This recovers the result (\ref{marutsutaresult}) by Maruhashi-Tsutaya.
\end{rem}

\section{Computing first cocycles}\label{compfirstcocycles}
We also solve a linear equation.

\subsection{Trivial representation}
We get the following lemma directly.
\begin{lem}
$H^1(\mathcal{F}_P)|_{I}=\mathbb{C}\omega_1$.
\end{lem}

\subsection{Corresponding to the lowest weight 1}
We also characterize the coboundary space. Before that, we prove that the special 1-cocycles are trivial.
\begin{lem}\label{oneside_1cocycle_trivial}
If $f_1\,\omega_1 \in Z^1(\mathcal{F}_P)|_{U^{-1}}$, then $f_1=0$.
\end{lem}
\begin{proof}
Put $g\, \omega_1\wedge\omega_2 = d_{\mathcal{F}_P}(f_1\,\omega_1)$ and $f_{1\, 0}=0$. Using (\ref{compdFP12_changed}), for any positive integer $N\geq 3$, we compute
\begin{align}
&\, \sum_{m=1}^{N-1}\sqrt{m}g_m \\
&= \sum_{m=1}^{N-1}\sqrt{m}\left( f_{1\,m-1}\sqrt{(m-1)m}-2mf_{1\,m}+f_{1\,m+1}\sqrt{m(m+1)} \right) \nonumber \\
&= \sum_{m=0}^{N-2}\sqrt{m+1}f_{1\, m}\sqrt{m(m+1)} + \sum_{m=1}^{N-1}\sqrt{m}(-2mf_{1\,m}) + \sum_{m=2}^{N}\sqrt{m-1}f_{1\, m}\sqrt{(m-1)m} \nonumber \\
&= \sum_{m=0}^{N-2}(m+1)\sqrt{m}f_{1\, m} + \sum_{m=1}^{N-1}(-2m)\sqrt{m}f_{1\,m} + \sum_{m=2}^{N}(m-1)\sqrt{m}f_{1\, m} \nonumber \\
&= 0 + \sum_{m=N-1}^{N-1}(-2m)\sqrt{m}f_{1\,m} + \sum_{m=N-1}^{N}(m-1)\sqrt{m}f_{1\, m} \nonumber \\
&= -N\sqrt{N-1}f_{1\, N-1} + (N-1)\sqrt{N}f_{1\, N}. \nonumber
\end{align}
Since $g=0$, it means
\begin{align}
f_{1\, N} = \sqrt{\frac{N}{N-1}}f_{1\, N-1}.
\end{align}
Thus we obtain
\begin{align}
f_{1\, N} = \sqrt{N}f_{1\, 1}.
\end{align}
Then $f_1=0$ because $\sum_{N=1}^{\infty}|f_{1\, N}|^2 <\infty$.
\end{proof}
\begin{lem}\label{1complem_1cocycle}
The following holds.
\begin{align}
B^1(\mathcal{F}_P)|_{U^{-1}}=\left\{f_1\,\omega_1+f_2\,\omega_2 \in Z^1(\mathcal{F}_P)|_{U^{-1}} \middle|\sum_{m=1}^{\infty}\sqrt{m}f_{2\,m}=0\right\}.
\end{align}
Specially, $Z^1(\mathcal{F}_P)|_{U^{-1}}$ is spanned by $\phi_1(\omega_1 - \omega_2)$ and $B^1(\mathcal{F}_P)|_{U^{-1}}$.
\end{lem}
\begin{proof}
Proving that the right hand side contains $B^1(\mathcal{F}_P)|_{U^{-1}}$ is easy. To prove opposite, let $N$ be a positive integer. We put
\begin{align}
h_{m}= \frac{N+1-m}{\sqrt{m(N+1)}}\,\,\,(1\leq m\leq N),
\end{align}
\begin{align}
h=\sum_{m=1}^{N}h_{m}\phi_m.
\end{align}
Then, we have
\begin{align}
d_{\mathcal{F}_P}h= (\mathrm{some\,\,function})\,\omega_1 + (-\sqrt{N+1}\phi_1+\phi_{N+1})\, \omega_2
\end{align}
by (\ref{compdFP01_f2_changed}). This formula and Lemma \ref{oneside_1cocycle_trivial} prove the result.
\end{proof}
As in Section \ref{1comp_2cocycle}, the following Sobolev estimation holds:
\begin{align}
||\xi||_k^2 \leq ||\eta||_{k+3}^2,
\end{align}
where $\eta=f_1\,\omega_1 + f_2\,\omega_2$ is an element of the right hand side in Lemma \ref{1complem_1cocycle} and $\xi=h$ is a $0$-cochain constructed as Lemma \ref{1complem}.

\subsection{Corresponding to the highest weight -1}
Similar argument also holds in this case.
\begin{lem}\label{-1complem_1cocycle}
The following holds.
\begin{align}
B^1(\mathcal{F}_P)|_{U^{1}}=\left\{f_1\,\omega_1+f_2\,\omega_2 \in Z^1(\mathcal{F}_P)|_{U^{1}} \middle|\sum_{m=-\infty}^{-1}(-1)^{-m}\sqrt{-m}f_{2\,m}=0\right\}.
\end{align}
Specially, $Z^1(\mathcal{F}_P)|_{U^{1}}$ is spanned by $\phi_{-1}(\omega_1 + \omega_2)$ and $B^1(\mathcal{F}_P)|_{U^{1}}$.
\end{lem}
As in Section \ref{-1comp_2cocycle}, the following Sobolev estimation holds:
\begin{align}
||\xi||_k^2 \leq ||\eta||_{k+3}^2,
\end{align}
where $\eta=f_1\,\omega_1 + f_2\,\omega_2$ is an element of the right hand side in Lemma \ref{-1complem_1cocycle} and $\xi=h$ is some $0$-cochain.

\subsection{The other cases}
Set $\pi\neq I, U^{-1}, U^1$ and fix $m_\pi \in\mathbb{M}$. We still start with proving the trivialness of the special 1-cocycles.
\begin{lem}\label{mod4_1cocycle_trivial}
Let $f_1\,\omega_1+f_2\,\omega_2 \in Z^1(\mathcal{F}_P)|_{\pi}$. Assume $f_{1\, m}=f_{2\,m}=f_{1\, m+3}=f_{2\,m+3}=0$ for any $m\in\mathbb{M}$ which satisfies $m\equiv m_\pi \,\,(\mathrm{mod}\,4)$. Then, $f_1=f_2=0$.
\end{lem}
\begin{proof}
Under the assumption, 1-cocycle conditions (\ref{compdFP12_changed}) is realized as the kernel of the linear map
\begin{align}
(f_{1\,m+1},f_{2\,m+1},f_{1\,m+2},f_{2\,m+2})\mapsto(g_{m},g_{m+1},g_{m+2},g_{m+3})
\end{align}
for each $m\equiv m_\pi$. (See (\ref{d12_matrix}).) Recall that the determinant $\gamma_{q\,m}$ is not $0$. Thus $f_1=f_2=0$.
\end{proof}
We consider the linear map
\begin{align}\label{44_zero_one_map}
(h_{m},h_{m+1},h_{m+2},h_{m+3})\mapsto(f_{1\,m+1},f_{2\,m+1},f_{1\,m+2},f_{2\,m+2})
\end{align}
for each $m\equiv m_\pi +2$ in (\ref{compdFP01_f1_changed}) and (\ref{compdFP01_f2_changed}). The coefficient matrix is also a block diagonal matrix whose block is a $4 \times 4$ matrix. Each block is represented as
\begin{align}\label{d01_matrix}
\begin{pmatrix}
f_{1\,m+1} \\
f_{2\,m+1} \\
f_{1\,m+2} \\
f_{2\,m+2}
\end{pmatrix}
=
\begin{pmatrix}
-\alpha_{q\,m} & 0 & \beta_{q\,m+2} & 0 \\
\alpha_{q\, m} & -2(m+1) & \beta_{q\,m+2} & 0 \\
0 & -\alpha_{q\,m+1} & 0 & \beta_{q\, m+3} \\
0 & \alpha_{q\,m+1} & -2(m+2) & \beta_{q\,m+3}
\end{pmatrix}
\begin{pmatrix}
h_{m} \\
h_{m+1} \\
h_{m+2} \\
h_{m+3}
\end{pmatrix}
.
\end{align}
Its determinant is also $\gamma_{q\,m}$ defined in (\ref{our_each_gamma}).
\begin{lem}
$H^1(\mathcal{F}_P)|_\pi =\{0\}$.
\end{lem}
\begin{proof}
Take any $f_1\,\omega_1+f_2\,\omega_2 \in Z^1(\mathcal{F}_P)|_{\pi}$. From (\ref{44_zero_one_map}), we can construct $h$ such that $f_1\,\omega_1+f_2\,\omega_2 - d_{\mathcal{F}_P}h$ satisfies the assumption of Lemma \ref{mod4_1cocycle_trivial}. Then $f_1\,\omega_1+f_2\,\omega_2 = d_{\mathcal{F}_P}h$.
\end{proof}
As in Section \ref{othercomp_2cocycle}, the following Sobolev estimation holds:
\begin{align}
||\xi||_k^2 \leq \Bigl( \frac{3!}{q_\Gamma^2}\Bigr)^2 2^{5(k+3)+8} ||\eta||_{k+3}^2,
\end{align}
where $\eta=f_1\,\omega_1 + f_2\,\omega_2$ is any $2$-cocycle and $\xi=h$ is some $0$-cochain.

\subsection{The whole sum}
As in Section \ref{whole_sum_2coboundaries}, for each $\eta\in B^1(\mathcal{F}_P)$, we have $\xi\in\Omega^0(\mathcal{F}_P)$ such that
\begin{align}
|| \xi ||_k^2 \leq C_{\Gamma \, k+3} || \eta ||_{k+3}^2.
\end{align}
The following holds.
\begin{prop}\label{1-cohomology}
The set $\{x,\,y_1,\, ...\, ,\, y_{2g}\}$ is basis for $H^1(\mathcal{F}_P)$, where the number $g$ is the multiplicity of $U^{-1}$ and $U^{1}$.
\end{prop}
\begin{rem}
This recovers the result (\ref{matsumitsuresult}) by Matsumoto-Mitsumatsu.
\end{rem}

\section{Determining the ring structure}\label{determinering}
We can prove our main theorem by combining the above preparation with the following lemma.
\begin{lem}\label{determineringlem}
Let $\phi_{1}, \phi_{1}'\in L^2(M_\Gamma)$ be weight vectors of $U^{-1}$. Here, $\phi_{1}$ and $\phi_{1}'$ do not necessarily belong to the same irreducible component. Also, let $\phi_{-1}, \phi_{-1}'\in L^2(M_\Gamma)$ be weight vectors of $U^{1}$. Then,
\begin{align}\label{ring_casimir_1}
(X_0^2 -X_1^2 -Y^2)(\phi_{1}\phi_{1}')=-2\phi_{1}\phi_{1}',
\end{align}
\begin{align}\label{ring_casimir_-1}
(X_0^2 -X_1^2 -Y^2)(\phi_{-1}\phi_{-1}')=-2\phi_{-1}\phi_{-1}',
\end{align}
\begin{align}\label{ring_x0}
X_0(\phi_{1}\phi_{-1})=0.
\end{align}
Especially, $\phi_{1}\phi_{1}'$, $\phi_{-1}\phi_{-1}'$ and $\phi_{1}\phi_{-1}$ orthogonal to $U^{-1}$ and $U^{1}$.
\end{lem}
\begin{proof}
Formulae are proved easily. The first two of them mean that $\phi_{1}\phi_{1}'$ and $\phi_{-1}\phi_{-1}'$ are eigenvectors corresponding to $-2$ of the Casimir element. On the other hand, the Casimir element vanishes on $U^{-1}$ and $U^1$. Then, they orthogonal to $U^{-1}$ and $U^{1}$. Also $\phi_{1}\phi_{-1}$ does by (\ref{ring_x0}). Indeed, the set $\mathbb{M}$ of $U^{-1}$ and $U^{1}$ does not contain $0$.
\end{proof}
\begin{proof}[Proof of Theorem \ref{MyThm2v3}]
Generators of $H^{*}(\mathcal{F}_P)$ are given in Proposition \ref{2-cohomology} and \ref{1-cohomology}. The vanishing of $y_i\wedge y_j$ in $H^2(\mathcal{F}_P)$ follows from Lemma \ref{determineringlem} for each $1\leq i,j \leq 2g$.
\end{proof}
\begin{rem}\label{M_Gamma_ring}
When $\Gamma$ is the fundamental group of a closed orientable hyperbolic surface, the vanishing of $y_i\wedge y_j$ is implied by the ring structure of $H^*_{\mathrm{dR}}(M_{\Gamma})$. In fact, we get the embedding $H^*_{\mathrm{dR}}(M_{\Gamma})/H^3_{\mathrm{dR}}(M_{\Gamma})\subset H^{*}(\mathcal{F}_P)$ as rings from (\ref{matsumitsuresult}) and (\ref{marutsutaresult}). The ring structure of $H^*_{\mathrm{dR}}(M_{\Gamma})$ is determined by the Thom-Gysin sequence and \cite[Lemma 1]{massey_ring}.
\end{rem}

\end{document}